\documentclass[12pt]{article}
\begin{document}
\title{OPTIMAL ACOUSTIC MEASUREMENTS}
\author{Margaret Cheney 
\thanks{Department of Mathematical Sciences,
Rensselaer Polytechnic Institute, Troy, NY 12180}
\and  David Isaacson 
\thanks{Department of Mathematical Sciences,
Rensselaer Polytechnic Institute, Troy, NY 12180}
\and Matti Lassas
\thanks{ Rolf Nevanlinna Institute, P.O. Box 4,
00014 University of Helsinki, FINLAND}}
\date{\today}
\maketitle

\begin{abstract}
We consider the problem of obtaining information about an inaccessible
half-space from acoustic measurements made in the accessible half-space.  If
the measurements are of limited precision, some scatterers will be undetectable
because their scattered fields are below the precision of the measuring 
instrument.  How can we make measurements that are optimal for detecting
the presence of an object?   In other words, what
incident fields should we apply that will result in the biggest measurements?  

There are many ways to formulate this question, depending on the measuring
instruments.  In this paper we consider a formulation involving 
wave-splitting in the accessible half-space:  what downgoing wave will
result in an upgoing wave of greatest energy?  

A closely related question arises in the case when we have a guess about the
configuration of the inaccessible half-space.  What measurements should we
make to determine whether our guess is accurate?   In this case we compare the
scattered field to the field computed from the guessed configuration.  
Again we look for the incident field that results in the greatest energy 
difference.

We show that the optimal incident field can be found by an iterative
process involving time reversal ``mirrors''.  For band-limited incident
fields and compactly supported scatterers, in the generic case
this iterative process converges to a single time-harmonic field.  
In particular, the process
automatically ``tunes" to the best frequency.  This analysis provides a
theoretical foundation for the frequency-shifting and 
pulse-broadening observed in certain
computations \cite{CT} and time-reversal experiments \cite{PF}
\cite{PTF}.

\end{abstract}

\section{Introduction}

This paper is motivated by the question ``What is the best way to do
acoustic imaging?''  If we want to make the best possible images, we
must begin with data that contain the most possible information.   
In particular, since all practical measurements are of limited precision,
some scatterers may be undetectable because their scattered fields are
below the precision of the measuring instrument:  our data will contain
no information about them.   What incident fields are ``best'', in
the sense that their scattered fields give the biggest measurements?  

This paper considers only the problem of detecting the presence of
an object (or distinguishing it from a guess) and not 
the problem of making an image of that object.  
For imaging, there are other criteria for ``best'' that one could
imagine using. 
  A Bayesian criterion \cite{L} \cite{LPS}, for example, 
would be to look for the measurement
  producing the ``narrowest'' posterior distribution for the scatterer
  when an {\it a priori} distribution for the scatterer is given.

The detection or distinguishability problem has been studied for fixed-frequency
problems in electrical impedance tomography \cite{I} and acoustic
scattering \cite{MNW}.  The connection between optimal measurements and
iterative time-reversal experiments was pointed out in \cite{MNW},
\cite{PF}, and  \cite{PTF}; in all these papers, the analysis was
carried out at a single fixed frequency.  The issue of optimal
time-dependent waveforms in a special $1+1$ -- dimensional case was
studied in  \cite{CT}, where a time-harmonic waveform was found to be
optimal.

In this paper we study the question of optimal time-dependent waveforms
in the $3+1$ -- dimensional case.  In particular, we consider the
half-space geometry:  we imagine that a plane divides space into
accessible and inaccessible regions, and we assume that we can make
measurements everywhere on the plane.

Section 2 contains a careful formulation of the idealized problem:  the
wave equation, the measurements, the notion of ``biggest''.  Section 3
is devoted to the example of a one-dimensional medium, in which the
problem can be solved explicitly.  Section 4 gives an iterative
experimental method that can be used to find the optimal field even if
the scatterer is unknown.  This method is precisely the iterative
time-reversal procedure of \cite{PF} and \cite{PTF}.  Section 5
discusses implications and open questions.  The paper concludes with
three appendices containing the technical details needed for the proof
of convergence of the time-reversal iterates.  We show that in general
the iterates converge to a time-harmonic field that is ``tuned" to the
best frequency.

\section{Basic Concepts}

\subsection{Distinguishability}
For any two operators $A_1$ and $A_2$, we say that $A_1$ is
distinguishable from $A_2$ with measurement precision $\epsilon$ if the
distinguishability $\delta(A_1, A_2)$, defined as 
\begin{equation}
\delta(A_1, A_2) = \sup_f {|| A_1 f - A_2 f || \over || f ||},
\label{distinguish} 
\end{equation} 
is greater than $\epsilon$.  A field that is best for distinguishing
$A_1$ from $A_2$ is an $f$ for which the maximum is attained.
We will determine below the norms that are
appropriate to use.

\subsection{Acoustic Wave Equation}

We consider the constant-density acoustic wave equation
\begin{equation}
(\nabla^2 - c^{-2}(x) \partial_t^2) U(t,x) = 0.\label{wave eq.}
\end{equation}
in the case in which $c=c_0$ everywhere in the upper half-space $x_3 > 0$.
This model includes neither dispersion nor dissipation.

We can formulate the scattering problem in a variety of ways \cite{CI}.
In particular, we can use either a boundary map, sources, or a scattering
operator defined in terms of wave splitting. 

\subsection{The Boundary Map}

To define the boundary map, we specify that $U=f$ on the surface
$x_3=0$.  This condition, together with an outgoing radiation 
condition at infinity \cite{CI},
uniquely determines a solution $U$ in the lower
half-space.  We can then take the normal derivative $\partial U /
\partial x_3$; this normal derivative, restricted to the surface
$x_3=0$, we denote by $g$.  The mapping from $f$ to $g$ is the boundary
map $\Lambda$.  Thus, on the surface $x_3 =0$,  $\Lambda U = \partial U
/ \partial x_3$.  Note that $\Lambda$ is an operator-valued function of
time.

Acoustic distinguishability can be defined in terms of the boundary map
as
\begin{equation}
\delta_B(c, c^0) = \sup_f {\| (\Lambda - \Lambda_0) f \| \over \|f\|}
\end{equation}
for appropriate norms.  Here $\Lambda_0$ denotes the boundary map for the
reference sound speed $c^0(x)$.  
This formulation, in terms of the boundary map,
is not pursued in this paper.

\subsection{Sources}

To formulate scattering in terms of sources, we consider 
the wave equation with a source:
\begin{equation}
(\nabla^2 - c^{-2}(x) \partial_t^2) U_J(t,x) = J(t, x).
\end{equation}
Scattering data is then $U_J(t,x)$ for $x$ on the plane, where $J$ is supported on
or above the plane.  
Acoustic distinguishability in terms of sources would be
\begin{equation}
\delta_S(c,c^0) = \sup_J { \| U_J - U_J^0 \| \over \| J \|},
\end{equation}
where $U_J^0$ denotes the field due to the reference sound speed
$c^0(x)$ and source $J$.  This formulation is not pursued in this
paper; instead we consider the scattering operator.

\subsection{The Scattering Operator}

We define the scattering operator in terms of upgoing and downgoing waves.  
The motivation
for this point of view is the existence of network analyzers, which can
decompose a time-harmonic signal in a waveguide into an upgoing one and a 
downgoing one, and measure
the amplitude and phase of the upgoing wave.  Stepped-frequency radar,
for example, is based on the ability of such instruments to transmit and 
receive signals at the same time.

\subsubsection{Upgoing and downgoing waves}

To define upgoing and downgoing waves, we make use of 
two Fourier transforms, a temporal one and a spatial one.
First we inverse-Fourier transform the solution $U$ of
(\ref{wave eq.}) in $t$: 
\begin{equation}
u(\omega, x) = {\mathcal F^{-1}} U = (2 \pi)^{-1} \int U(t,x) e^{i \omega t} dt.
	\label{tFT}
\end{equation}
This frequency-domain solution
$u$ satisfies the reduced wave equation
\begin{equation}
(\nabla^2 + \omega^2 c^{-2}) u(\omega,x) = 0. \label{fwave eq.}
\end{equation}
We write $k = \omega / c_0$ and with a small abuse of notation
we write $u(k,x)$ instead of $u(c_0 k, x)$.
We then Fourier transform $u$ again in $x' = (x_1, x_2)$, so that
\begin{equation}
\hat u(k, \eta', x_3) = F_{x_3} u = \int u(k, x) 
	e^{-i k \eta' \cdot x'} d^2 x' \label{spaceFT}
\end{equation}
where $\eta' = (\eta_1, \eta_2)$.  (We note that $F_{x_3}$ depends
also on $k$.)  Then $U$ is recovered as
\begin{equation}
U(t,x) = {1 \over (2 \pi)^2} \int \int \hat u(k,\eta',x_3)
	 e^{ik\eta' \cdot x'} e^{-i k c_0 t} k^2 d^2 \eta' c_0 dk.
\label{FT formula}
\end{equation}

In the upper half-space, $\hat u$ satisfies the ordinary differential equation
\begin{equation}
(\partial_{x_3}^2 + k^2 - k^2 |\eta'|^2) \hat u = 0,
\end{equation} 
which has the general solution
\begin{equation}
\hat u (k,\eta',x_3) = A(k,\eta') e^{i k\eta_3 x_3} + 
	B(k,\eta') e^{-ik \eta_3 x_3}, \label{gen.sol.}
\end{equation}
where
\begin{equation}
\eta_3 (k) = \cases{ \sqrt{1- |\eta'|^2} & for $1 > |\eta'|$ \cr
            i  ({\rm sgn} k) \sqrt{|\eta'|^2 - 1} & for $1 < |\eta'|$ \cr}
	 \label{xi3}
\end{equation}
We define the vectors $\eta^\pm = (\eta', \pm \eta_3)$, which satisfy 
$\eta^\pm \cdot \eta^\pm = 1$.

In order for $U$, as defined by (\ref{FT formula}),
to be real-valued, the Fourier transform $\hat u$ must satisfy certain 
symmetry conditions.  In particular, we must have  
$A(-k,\eta') = \overline{ A(k,\eta')}$,
$B(-k,\eta') =\overline{ B(k,\eta')}$, and
$\eta_3 (-k) = \overline{\eta_3 (k)}$.

Equation (\ref{gen.sol.}) shows us how to split 
the time-domain solution of (\ref{wave eq.}) into two parts, which
we call the upgoing and downgoing parts.  Thus for 
$x_3 > 0$ we write
$U = U^\uparrow + U^\downarrow$, where $U^\downarrow$ is
\begin{equation}
U^\uparrow (t,x) = \int \int A(k, \eta' ) e^{i k \eta^+ \cdot x}
	e^{-i k c_0 t} k^2 d^2 \eta' c_0 dk,  \label{up}
\end{equation} 
and
\begin{equation}
U^\downarrow (t,x) = \int \int B(k, \eta' )  e^{i k \eta^- \cdot x} 
	e^{-i k c_0  t} k^2 d^2 \eta'  c_0 dk. \label{down}
\end{equation}

We see that (\ref{down}) and (\ref{up}) are plane wave decompositions.  
The components for which $|\eta'|<1$ are propagating plane waves, and $\eta^\pm$ 
is a unit vector that gives the direction of propagation.  The sign of the third
component of $\eta^\pm$ determines whether the wave is downgoing or upgoing.  
On the other hand, components with $|\eta'| > 1$ correspond to evanescent
waves.  For $U^\downarrow$, these evanescent waves
decay in the downward (negative $x_3$) direction; 
for $U^\uparrow$, they decay in the upward (positive) direction.

\subsubsection{The scattering operator}

It is natural to define a scattering operator ${\mathcal S}^\updownarrow$ 
as the map from $U^\downarrow$ to $U^\uparrow$.  We denote the kernel of
this operator also by ${\mathcal S}^\updownarrow$:
\begin{equation}
U^\uparrow (x,t) = 
\int_{-\infty}^\infty \int {\mathcal S}^\updownarrow(x,y, t-\tau)  U^\downarrow (y, \tau) d^3y d\tau.
\label{St}
\end{equation}
We note that this scattering operator is defined only on downgoing solutions
of the Helmholtz equation, i.e., on functions of the form (\ref{down}).

The kernel of (\ref{St}) is a convolution in time because the 
Fourier transform ``diagonalizes"
the time derivative of (\ref{wave eq.}), so that the frequency 
is simply a parameter in (\ref{fwave eq.}).  
In other words, the convolution is an expression of the fact
that in the frequency domain, (\ref{St}) takes the form  
\begin{equation}
u^\uparrow (k, x) = \int S^\updownarrow(k, x, y)
	 u^\downarrow (k, y) d^3y.
	\label{Somega}
\end{equation}
The time-domain operator ${\mathcal S}^\updownarrow$  is related to
the frequency-domain scattering operator $S^\updownarrow$ by 
$S^\updownarrow = {\mathcal F^{-1} \mathcal S^\updownarrow \mathcal F}$.

In fact, $\mathcal S^\updownarrow$ and $S^\updownarrow$ are determined by their
actions on the plane $x_3 = 0$.  We see this as follows.
First we Fourier transform (\ref{Somega})  in space.  The 
operator $S^\updownarrow$ is transformed into the operator  
$\hat S_{x_3, \tilde x_3} = F_{x_3}S^\updownarrow F^{-1}_{\tilde x_3} = 
(F_{x_3} {\mathcal F}^{-1}) {\mathcal S}^\updownarrow 
(F_{\tilde x_3}{\mathcal F}^{-1})^{-1}$,
where $F_{x_3}$ is defined by (\ref{spaceFT}) and 
${\mathcal F}^{-1}$ by (\ref{tFT}).  
The transformed version of (\ref{Somega}) is 
\begin{equation}
\hat u^\uparrow (k, \eta', x_3) = 
	\int \hat S_{x_3, \tilde x_3} (k, \eta', \tilde \eta') 
	\hat u^\downarrow (k, \tilde \eta', \tilde x_3)
	k^2 d^2 \tilde \eta'. \label{FTS}
\end{equation}
The operator $\hat S_{x_3, \tilde x_3}$ is determined
completely by its action at $x_3=0$, which we see from the following argument.

Into (\ref{FTS}) we substitute
$\hat u^\downarrow(k, \tilde \eta', \tilde x_3) 
	= B(k, \tilde \eta') e^{-i k \tilde \eta_3 \tilde x_3}$ 
and 
$\hat u^\uparrow(k, \eta', x_3)
	 = A(k, \eta') e^{ik \eta_3 x_3}$; 
we see that 
\begin{equation}
\hat  S(k, \eta', \tilde \eta')  =  
	e^{-ik \eta_3 x_3} 
	\hat S_{x_3, \tilde x_3} (k, \eta', \tilde \eta') 
	 e^{-ik \tilde \eta_3 \tilde x_3} \label{hat S}
\end{equation}
satisfies
\begin{equation}
A(k, \eta') = \int \hat S(k, \eta', \tilde \eta') 	
	B(k, \tilde \eta') k^2 d^2 \tilde \eta'
\label{hat s}
\end{equation}
where $A$ and $B$ are as in (\ref{gen.sol.}).
The relation between the operators $\hat S$ and 
$\hat S_{x_3, \tilde x_3}$ can be written
$\hat S_{x_3, \tilde x_3} = E_{x_3} \hat S  E_{\tilde x_3}$,
where $E_{x_3} $ is the operator of multiplication by $\exp ( ik \eta_3 x_3)$,
and thus 
\begin{equation}
S^\updownarrow = F_{x_3}^{-1} E_{x_3} \hat S E_{\tilde x_3} F_{\tilde x_3}.
\end{equation}
On the plane $x_3=0$, this becomes 
\begin{equation}
S = F_0^{-1} \hat S F_0 \label{S}
\end{equation}
This defines a scattering operator $S$ on the plane $x_3 = 0$. 
 It is this operator,
together with the corresponding time-domain operator 
$\mathcal S = \mathcal F S \mathcal F^{-1}$,
that we will use in the rest of the paper.  	
We note that the domain of the operator
$\mathcal S$ is restricted to the space of downgoing waves as defined by
 (\ref{down}).

The scattering operator $\hat S$ and the boundary map $\Lambda$ 
are related to each other by formulas developed in \cite{CI} 
(See Appendix A for details).  
They are thus equivalent operators.  
Whether it is better to formulate a given problem
in terms of a scattering operator or a boundary map depends largely on the
design of the equipment involved.

\subsection{The Energy Identity and the Energy Flux}

If we multiply (\ref{wave eq.}) by $\partial_t U$ and integrate
the resulting equation over the volume $V$, we obtain
\begin{equation}
\int_V \bigg( (\partial_t U) \nabla^2 U - 
  {1 \over 2 c^2(x)} \partial_t (\partial_t U)^2 \bigg) dx = 0.
\label{Ut*wave}
\end{equation}
We write the first term of (\ref{Ut*wave}) as 
$\nabla \cdot ((\partial_t U) \nabla U) -
 \nabla (\partial_t U) \cdot \nabla U$,
and apply the divergence theorem to the term containing the divergence.
We thus obtain
\begin{equation}
\int_{\partial V} (\partial_t U) \partial_\nu U dS =  \partial_t 
  \int_V {1 \over 2} \bigg( |\nabla U|^2 + {1 \over c^2(x)}
 (\partial_t U)^2 \bigg) dx ,
\label{energy identity}
\end{equation}
where $\nu$ denotes the outward unit normal to the surface $\partial V$.

This equation relates the change in energy in the volume $V$
(the right side of (\ref{energy identity})) to the energy flux across
its boundary surface $\partial V$.  

From (\ref{energy identity}) we see that the time-integrated energy flux
 across a
surface $\partial V$ in the normal direction $\nu$ is
\begin{equation}
W(U) = 
 - \int_{-\infty}^{\infty} \int_{\partial V} (\partial_t U)
        \partial_{\nu} U dS dt .
\label{flux}
\end{equation}
We can use Parseval's identity to write the time-integrated energy flux
in terms of the frequency-domain wave functions: 
\begin{equation}
W(u) =  -(2 \pi)^3 \int_{-\infty}^{\infty} \int_{\partial V} (i \omega u)
        \overline{ \partial_{\nu} u }\ dS d\omega
	= -(2 \pi)^3 \int_{-\infty}^{\infty} \int_{\partial V} (i c_0 k u)
        \overline{ \partial_{\nu} u }\ dS c_0 dk, \label{Parseval}
\end{equation}
where the overline denotes the complex conjugate.

\subsubsection{The energy flux of upgoing and downgoing waves}

At the surface $x_3 = 0$, the total field splits into upgoing and downgoing 
parts.  Because the flux is quadratic, it does not necessarily split
into corresponding upgoing and downgoing fluxes.  However, a quick 
calculation using (\ref{up}), (\ref{down}), and 
(\ref{Parseval}) shows that if upgoing and downgoing evanescent waves
are not both present on the plane $x_3=0$, the time-integrated cross terms 
$ \int_{-\infty}^{\infty} \int_{x_3 = 0} (\partial_t U^\downarrow)
        \partial_{x_3} U^\uparrow dx' dt$ and
$ \int_{-\infty}^{\infty} \int_{x_3 = 0} (\partial_t U^\uparrow)
        \partial_{x_3} U^\downarrow dx' dt$ cancel.
Under these conditions, the time-integrated fluxes do split into upgoing
and downgoing fluxes, so that 
$W(U^\downarrow+U^\uparrow) = W(U^\downarrow)+W(U^\uparrow)$.  
Throughout this paper we assume that the sources of the downgoing
field are far from the scatterers, so that there is no interaction
between upgoing and downgoing evanescent waves on the plane $x_3=0$.
The flux of a downgoing wave is positive; that of an upgoing wave
is negative.

 We write the downgoing energy flux as 
\begin{equation}
W(U^\downarrow) = 
 \int_{-\infty}^{\infty} \int_{x_3 = 0} (\partial_t U^\downarrow)
        \partial_{x_3} U^\downarrow dx' dt
= (2 \pi)^3 \int_{-\infty}^{\infty} \int |B(k,\eta')|^2 c_0^2  \eta_3
 	k^4 d^2 \eta' dk , 
\label{Uinc}
\end{equation}
where we have used (\ref{down}) and (\ref{Parseval}) 
in carrying out the computation (\ref{Uinc}).
In (\ref{Uinc}) there is no minus sign because downgoing energy
travels in the $-x_3$ direction.

Although the left side of (\ref{Uinc}) is real, it is not obvious that the
right side is, because $\eta_3$ can be imaginary.  However, if one splits
the $k$ integral into pieces as
\begin{equation}
W(U^\downarrow) = \int \bigg( \int_{-\infty}^0 + \int_0^\infty \bigg)
	|B(k,\eta')|^2 c_0^2 k^4\eta_3 dk d^2 \eta'
\end{equation}
and uses the symmetry properties of $B$ and $\eta_3$, one sees that
for evanescent waves, the two terms cancel. 
 This shows that the evanescent waves do not
contribute to the energy flux.

The flux can be used to form an inner product on the space of 
downgoing propagating waves; we define
\begin{eqnarray}
(U^\downarrow, V^\downarrow)_W 
	&=& {1 \over 2} \int_{-\infty}^\infty  \int_{x_3 = 0} 
	\left( (\partial_t U^\downarrow) \partial_{x_3} V^\downarrow +
	(\partial_t V^\downarrow) \partial_{x_3} U^\downarrow \right) 
	d^2x' dt \cr
	&=& (2 \pi)^3 \int \int_{|\eta'|<1} \hat u^\downarrow 
	\overline{\hat v^\downarrow}
	c_0^2 k^4 \eta_3 d k d^2 \eta' \label{upfluxprod}
\end{eqnarray}
We note that the product $c_0^2 k^4 \eta_3$ is non-negative, so in the transform
domain, this inner product is merely a weighted $L^2$ inner product.
 
Similarly, the energy flux of the upgoing scattered field 
$S U^\downarrow = U^\uparrow$ that passes through the plane $x_3=0$ is
\begin{equation}
W(U^\uparrow)= \int_{-\infty}^{\infty} \int_{x_3 = 0} (\partial_t U^\uparrow)
        \partial_{x_3} U^\uparrow dx' dt
= - (2 \pi)^3 \int_{-\infty}^{\infty} \int |A(k,\eta')|^2 
	c_0^2 k^4 \eta_3 d^2 \eta' dk.
\label{Usc}
\end{equation}
The minus sign in (\ref{Usc}) is due to the fact that the upgoing
wave corresponds to energy leaving the lower half-space.  
Equations (\ref{Uinc}) and (\ref{Usc}) show that the time-integrated
downgoing flux $W(U^\downarrow)$ is positive and the time-integrated
upgoing flux $W(U^\uparrow)$ is negative. 

Note that for propagating waves, $|W|$ satisifes the triangle inequality:  
$|W(U^\uparrow_1 + U^\uparrow_2)| \leq |W(U^\uparrow_1)| + |W(U^\uparrow_2)|$.  

The flux inner product on the space of upgoing propagating waves is
\begin{eqnarray}
(U^\uparrow, V^\uparrow)_W 
	&=& \left| {-1 \over 2} \int_{-\infty}^\infty  \int_{x_3 = 0} 
	\left( (\partial_t U^\uparrow) \partial_{x_3} V^\uparrow +
	(\partial_t V^\uparrow) \partial_{x_3} U^\uparrow \right) \right| 
	d^2x' dt\cr
	&=&  (2 \pi)^3 \int \int_{|\eta'|<1} \hat u^\uparrow 
	\overline{\hat v^\uparrow}
	c_0^2 k^4 \eta_3 d k d^2 \eta';  \label{fluxprod}
\end{eqnarray}
thus for propagating waves, $|W(U^\uparrow)| = (U^\uparrow, U^\uparrow)_W$.  

\paragraph{Conservation of energy.}
If the medium is initially
quiescent, conservation of energy tells us that $W(U^\downarrow) \geq
|W(\mathcal S U^\downarrow)|$; this can be seen from integrating (\ref{energy identity})
over all time, and using the fact that the energy within the volume $V$
is initially zero, and cannot become negative.  The time integral of the
right side of (\ref{energy identity}) is thus positive.  The left side we
write as $W(U) = W(U^\downarrow) + W(\mathcal S U^\downarrow)$.  This implies that
the total upgoing flux $|W(\mathcal S U^\downarrow)|$ cannot be greater than 
the total downgoing flux $W(U^\downarrow)$.  

\paragraph{Finite-energy fields on the plane.}
We define the space $w$ of finite-energy functions on the plane to be the
closure of $C^\infty_0({\bf R}^2\times{\bf R})$ in the inner product
\begin{equation}
(u,v)_w^2= \int \int \hat u(k, \eta)
   \overline{\hat v(k, \eta)} c_0^2 k^4|\eta_3| d^2\eta'dk,
	\label{Wspace}
\end{equation}
and the space $W = {\mathcal F} w$.

\subsection{Acoustic Distinguishability via the Scattering Operator}

We define the acoustic distinguishability in terms of the upgoing and
downgoing energy fluxes through the surface $x_3=0$. 

For a reference scatterer with scattering operator ${\mathcal S}_0$,
the energy flux of the upgoing scattered field ${\mathcal S}_0
U^\downarrow$ and of the difference field $({\mathcal S}-{\mathcal
S}_0) U^\downarrow$ are defined similarly.

In general the distinguishability of ${\mathcal S}$ from ${\mathcal S}_0$ 
with the incident field $U^\downarrow$ is
\begin{equation}
\delta({\mathcal S, \mathcal S}_0) = \sup_{U^\downarrow}
	{|W(({\mathcal S} -{\mathcal S}_0) U^\downarrow)| \over W(U^\downarrow)}
	=\sup_{u^\downarrow} {|W((S-S_0)u^\downarrow)| \over W(u^\downarrow)} .
\label{Sdisting}
\end{equation}

We recall that evanescent components do not contribute to the energy flux.  To
remove the evanescent components from (\ref{Sdisting}), we denote by $P$ the 
orthogonal projection onto the propagating components:  
$P = F_0^{-1} \hat P F_0$, where $\hat P$ is the operator of multiplication by 
$\chi_{|\eta'|<1}$, the function that is 
one for $|\eta'|<1$ and zero otherwise.  Explicitly, $P$ is given by
\begin{equation}
P(k)f(x') = \int k^2 \int_{|\eta'|<1} e^{ik \eta' \cdot (x' - y')}
	 d^2 \eta' f(y') d^2 y'.
\end{equation}
In the time domain, $\mathcal P = \mathcal F P \mathcal F^{-1}$.
With this notation, we can write $W(U) = W(\mathcal P U) = W(Pu)$.  

Moreover, the scattered field due to an evanescent incident wave has zero total energy
flux.  This is because of the comments at the end of the
previous section:  $0 = W((I- P)u^\downarrow) \geq  |W(S(I- P) u^\downarrow)|$
implies that $W(S(I-P) u^\downarrow) = 0$.

In addition, the upgoing energy flux can only be increased by getting
rid of the evanescent components of the incident wave.  This is because
of the triangle inequality $|W(SPu + S(I-P)u)| \leq |W(SPu)| +
|W(S(I-P)u| = |W(SPu)|$.

This implies that the downgoing waves that give rise to the maximum total energy flux 
are propagating waves.  Thus we find that the distinguishability can be written
\begin{equation}
\delta(S, S_0) = \sup_{U^\downarrow} {|W(\mathcal P (\mathcal S - \mathcal S _0)
	\mathcal P U^\downarrow)| \over W(\mathcal P U^\downarrow)}
	= \sup_{u^\downarrow} {|W(P(S-S_0)P u^\downarrow)| \over
	W(Pu^\downarrow)}.  \label{PSdisting}
\end{equation}

We note that the scattering operator ${\mathcal S}_0$ for free space is the zero
operator.  Thus, according to (\ref{distinguish}) and 
(\ref{Sdisting}), the presence of a
scatterer can be detected with measurement precision $\epsilon$ if the
distinguishability satisfies 
\begin{equation}\delta({\mathcal S},0) = \sup_{U^\downarrow}
	{|W({\mathcal P \mathcal S \mathcal P} U^\downarrow)| 
	\over W(\mathcal P U^\downarrow)}
	= \sup_{\mathcal P U^\downarrow}{|W(\mathcal P U^\uparrow) |
	\over W({\mathcal P} U^\downarrow) } > \epsilon . \label{distingS}
\end{equation}

The distinguishability can be defined equally well in terms of the
operator $\hat S$ of (\ref{hat S}) or $S$ of (\ref{S}).

\medskip
\section{Example: The One-Dimensional Case}

If the medium in the lower half-space depends only on depth, then the 
coefficient $B$ of (\ref{gen.sol.})
is the reflection coefficient $R(k,\eta')$ multiplied by the 
incident coefficient $A$.  
In this case, the distinguishability $\delta({\mathcal S},0)$ can be computed from 
(\ref{Uinc}), (\ref{Usc}), and (\ref{distingS}) as
\begin{equation}
\delta({\mathcal S},0) 
	= \sup_{B}{\int \int |R(k,\eta') 
	B(k, \eta')|^2  c_0^2 k^4 \eta_3 dk d^2\eta'
	 \over \int \int |B(k,\eta')|^2 
	c_0^2 k^4 \eta_3 dk d^2\eta'}.
	\label{distingS2}
\end{equation}
The maximum of the right side of (\ref{distingS2}) is attained in the
limit when $B$ is a delta function supported at the maximum of $|R|$.

Thus to maximize the scattering from a one-dimensional scatterer, we
compute the conventional reflection coefficient $R(k, \eta')$, and find
the values of $k$ and $\eta'$ at which it attains its maximum.  Taking
$B$ to be a delta function supported at these points corresponds to
taking an incident field that is a plane wave of fixed frequency 
$\omega= c_0 k$ and incident direction given by $\eta'$.

Note that since $R$ is minus one for $|\eta'| = 1$ (grazing), a maximum
always occurs at grazing incidence.  If this is undesirable, grazing
incidence can be excluded by modifying the definition of
distinguishability.

\medskip

\section{An Adaptive Method for Producing the Best Fields}

To maximize the distinguishability when the medium is unknown, we can use
the following adaptive method.   

We write
\begin{equation}
\delta(\mathcal S, \mathcal S_0) = \sup_{U^\downarrow} 
	{|W(\mathcal P (\mathcal S - \mathcal S_0)
	\mathcal P U^\downarrow)| \over W(\mathcal P U^\downarrow)}
	= \sup_{U\in W} { ( U, 
	(\mathcal P (\mathcal S - \mathcal S_0) \mathcal P)^*
	\mathcal P (\mathcal S - \mathcal S_0) \mathcal P U)_W
\over
	 (\mathcal P U, \mathcal P U)_W}
\label{maxflux}
\end{equation}
where $U(x,t)=U^\downarrow|_{x_3=0}$ and where
the adjoint $^*$ has been taken in the space $W$ (defined just below
(\ref{Wspace})).  

We see in Appendix A that 
$(\mathcal P \mathcal S \mathcal P)^* = T (\mathcal P\mathcal S \mathcal P) T$,
 where $T$ denotes the time-reversal operator $TU(t,x) = U(-t,x)$.

Thus we see that the operator appearing on the right side of (\ref{maxflux}) is
$\mathcal A = T\mathcal P (\mathcal S - \mathcal S_0) \mathcal P T (\mathcal S - \mathcal S_0)$.  
In general, to maximize a quotient of the form 
\begin{equation}
\langle U, {\mathcal A}U \rangle / \langle U, U \rangle, \label{RRquotient}
\end{equation}
one considers an appropriately normalized sequence $\mathcal A^n U$.
When $\mathcal A$ is compact, this sequence converges to the largest
eigenvalue of $\mathcal A$.  Here, however, $\mathcal A$ has a
continuous spectrum, so we expect the sequence $\mathcal A^n U$, when
appropriately normalized, to converge to a generalized eigenfunction of
$\mathcal A$, and the corresponding quotient (\ref{RRquotient}) to
converge to the supremum of the continuous spectrum.  We note that
such generalized eigenfunctions do not have finite energy.

When $\mathcal A$ is compact, the usual way to normalize $\mathcal A^n
U$ is to divide by $\| \mathcal A^n U \|$.  Here, however, because we
expect ${\mathcal A}^n U$ to converge to a distribution in the time
variable, we must use a distributional normalization.  We consider test
functions in a particular space that is discussed in the appendix.
These test functions are functions of space and time.  In the time
variable, they are Fourier transforms of functions of compact
support.   The distribution action is chosen to coincide with the flux
inner product defined by (\ref{upfluxprod}) and (\ref{fluxprod}).  For
the distribution action we use the same notation $( \cdot, \cdot)_W$ as
for the flux inner product.

To normalize, we choose an arbitrary test function
$\Psi$, and consider the sequence 
$\mathcal A^n U /  (\mathcal A^n U, \Psi )_W$.
This gives rise to the following algorithm for carrying out the maximization of (\ref{maxflux}).
\medskip

1.  Start with any $V^\downarrow_0$; let $j=0$.

2.  Send $V^\downarrow_j$ into the lower half-space, and measure the resulting
upgoing field $V^\uparrow_j (t,x) = \mathcal S V^\downarrow_j (t,x)$.

3.  Calculate the corresponding scattering from the reference configuration
$\mathcal S_0 V_j^\downarrow (t,x))$.  Calculate the difference field
$\tilde V_j^\uparrow (t,x) = V^\uparrow_j(t,x) - \mathcal S_0 V_j^\downarrow (t,x)$.
    
4.   If $j$ is even, let 
\begin{equation}
V^\downarrow_{j+1}(t,x)  =  \tilde V_j^\uparrow (-t,x),
\end{equation}
add one to $j$, and go to Step 2.

5.  If $j$ is odd, normalize:
\begin{equation}
V^\downarrow_{j+1}(t,x)  = { \tilde V_j^\uparrow (-t,x) \over 
	(T \tilde V_j^\uparrow , \Psi )_W },
\end{equation}
add one to $j$, and go to Step 2.

\medskip

Appendix B contains a proof that, in the case of a compactly supported
scatterer in free space, the sequence 
$U_n = V_{2n} = \mathcal A^n U$ generally converges to a single
time-harmonic wave.  The frequency of this wave is the frequency at
which the largest eigenvalue of $S$ attains its maximum.  If this largest
eigenvalue happens to attain the same maximum at several different
frequencies, then the iterates $U_n$ converge to a sum of time-harmonic
waves with these frequencies.   The relative strengths of the different
frequencies is determined by the corresponding frequency components of
the initial incident wave $U_0 = V_0^\downarrow$. 

The argument in Appendix B takes place within a limited frequency band; this
frequency band is determined by the bandwidth of the test function.  

We note that as expected, the limiting time-harmonic waves do not have finite energy.  
This property also appears in the one-dimensional example (\ref{distingS2}).  

Step 4 can be omited and Step 5 performed for every $j$:  the linearity of 
the problem implies that extra normalizations do not affect the limit.  
The proof in the Appendix,
however, corresponds to the above algorithm.  

The algorithm can also be implemented including a step in which the evanescent
waves are filtered out.  If they are not filtered out, however, 
they will die out anyway as the iteration
proceeds, because experimental time-reversal of a field that 
includes evanescent 
waves is simply another physical field with evanescent waves.  

\section{Conclusions and Open Questions}

This analysis shows that the iterative time-reversal work of \cite{PF} 
and \cite{PTF} provides an
experimental method to obtain optimal fields.  Moreover, this analysis
explains the frequency-shifting and 
pulse-broadening seen in \cite{PTF} and \cite{CT}: the 
optimal time-domain waveform is a time-harmonic one tuned to the
best frequency.  

This analysis suggests that the commonly-used pings and chirps 
are not optimal from the point of view of distinguishability.

There are many open questions related to this work, 
one of which is the question of limited-aperture 
and limited-time measurements.
Upgoing and downgoing waves in a limited aperture
can be defined with the help of eigenvalues of the Laplacian for the
aperture.  However, it is not clear how to determine the entire
incident wave if the incident wave is known in only a limited aperture.  
This involves a detailed modeling of the transducer or antenna.
Perhaps a formulation in terms of sources will be more useful in this case.

We have not studied the question of whether the distinguishability,
as a function of the medium, is monotone in any sense.  This is an
important issue for the following reason.  Suppose we discover that a sphere
of a certain radius is detectable with a certain measurement
precision.  Does this imply that a larger object will also be
detectable?  For fixed-frequency measurements, the answer to this
question is certainly no, because of the phenomenon of resonance.  A
small sphere may happen to have a radius commensurate with the
wavelength of the probing wave, and may therefore scatter much more
strongly than a larger sphere.  However, the use of time-dependent fields 
may give different results.

The simple wave equation studied in this paper does not include the important
effects of variable density, dispersion, and dissipation.  

Moreover, the question of distinguishability is only the first step
in building an optimal imaging system.  How should we choose a full set of
optimal fields that could be used to form an image?

\section{Acknowledgments}
This work was partially supported by the Office of Naval Research.
M.C. would like to thank a number of people for helpful discussions:  
Gerhard Kristensson and his group in Lund,  Jim Rose, Claire
Prada, and Isom Herron;  M.L. thanks Lassi P\"aiv\"arinta for
interesting discussions.

\appendix
\section{Appendix:  Properties of $S$. }
\newtheorem{theorem}{Theorem}
\newtheorem{proposition}{Proposition}
\newtheorem{corollary}{Corollary}
\newtheorem{lemma}{Lemma}

\subsection{Expression for Kernel}

We can find an expression for the kernel 
$\hat S$ of (\ref{FTS}) and (\ref{hat S})
by taking $\hat u^\downarrow$ to be a delta function. 
The kernel of $\hat S$ is then the corresponding upgoing wave $\hat u^\uparrow$.
Taking $\hat u^\downarrow$ to be a delta function means that we take 
$u^\downarrow$ of the form $\exp (i k \tilde \eta^- \cdot x)$ 
for some $\tilde \eta^- = (\tilde \eta', -\tilde \eta_3)$.  
Here $\tilde \eta_3$ can be complex.  
We write the corresponding frequency-domain field as $\psi$:
\begin{equation}
\psi(k, x, \tilde \eta^-) = {1 \over 2 \pi} \int U(t,x, \tilde \eta^-) 
	e^{ik c_0 t} dt
\end{equation}
so that 
\begin{equation}
\hat u(k, \eta', x_3, \tilde \eta^-) = 
	\int \psi(k, x, \tilde \eta^-) e^{-ik \eta' \cdot x'} d^2 x'.
\end{equation}  

In the case of scattering from a perturbation in free space, 
we can express the 
total frequency-domain field $\psi$ as a solution of the Lippmann-Schwinger equation
\begin{equation}
\psi(k, x, \tilde \eta) = \exp (ik  \tilde \eta \cdot x) - 
	k^2 \int g(k,x,y) V(y) \psi(k, y, \tilde \eta) d^3 y, \label{LS}
\end{equation}
where $g$ is the usual outgoing Green's function
\begin{equation}
g(k,x,y) = {e^{ik |x-y|} \over 4 \pi |x-y|} \label{Green}
\end{equation}
and $V(y) = 1 - c_0^2/c^2(y)$.
The scattered field $u^\uparrow$ thus is represented by the integral term of (\ref{LS}).

The Green's function can be written in terms of its two-dimensional 
Fourier transform as \cite{CI}
\begin{equation}
g(k,x,y) =  {1 \over (2 \pi)^2} \int {i \over 2k \eta_3}
	 e^{ik \eta_3 |x_3 - y_3|} 
	e^{ik \eta' \cdot (x' - y')} k^2 d^2 \eta' .	\label{FG}
\end{equation}

To compute $A$ and $B$ of (\ref{gen.sol.}), we take the $x_1, x_2$ Fourier
transform of (\ref{LS}) in the region $x_3 > 0$. 
 We assume that the perturbation $V$ is supported in 
the region $y_3 < 0$, so that when we consider (\ref{LS}) we can remove 
the absolute values in (\ref{FG}).   
In the transform domain, the scattered field is given by 
\begin{equation}
\hat u^\uparrow (k, \eta', x_3, \tilde \eta^-)
	 = - k^2 \int{i \over 2k \eta_3}  
	\int e^{-ik \eta^+ \cdot y} V(y) \psi (k,y,\tilde \eta^-)
	 d^2 y' dy_3  e^{ik \eta_3 x_3}.
\end{equation}

This shows that 
\begin{equation}
\hat S(k, \eta', \tilde \eta') = -{ik \over 2 \eta_3} 
	A(k, \eta^+, \tilde \eta^-) \label{kernel}
\end{equation}
where 
\begin{equation}
A(k, \eta, \tilde \eta) = \int e^{-ik \eta \cdot y} V(y) 
	\psi (k, y, \tilde \eta) d^3 y
\end{equation}
is a scalar multiple of the classical scattering amplitude \cite{N}.  
(This $A$ is not to be confused with the $A$ of (\ref{gen.sol.})!)
It satisfies the reciprocity relation 
\begin{equation}
A(k, \eta, \tilde \eta) = A(k, -\tilde \eta, -\eta)
	\label{reciprocity}
\end{equation}
and the symmetry relation (for real-valued perturbations $V$)
\begin{equation}
\overline{ A(k, \eta, \tilde \eta)} = A(-k, \eta, \tilde \eta). 
	\label{symmetry}
\end{equation}

It is clear from (\ref{kernel}) that $\hat S$ is an analytic function
 of $k$ \cite{N}.

An expression similar to (\ref{kernel}) can be obtained for the field
scattered from a perturbed half-space or layered medium; in this case
the appropriate background Green's function should be used instead of
the free-space Green's function in (\ref{LS}).  Scattering theory in
such cases is considered, for example, in \cite{W}, \cite{X}, and
\cite{ER}.

\subsection{The Adjoint}

For a scatterer in free space, the adjoint of 
$\mathcal P \mathcal S \mathcal P$ 
in the space $W$ can be computed explicitly as follows.  
From (\ref{fluxprod}) we have
\begin{eqnarray}
( \mathcal P \mathcal S \mathcal P U, V )_W &=&
	(2 \pi)^3 \int_{-\infty}^\infty 
	\int_{|\eta'| <1} \int_{|\tilde \eta'|<1} 
	\hat S(k, \eta', \tilde \eta') 
	\hat u (k, \tilde \eta') k^2 d^2 \tilde \eta' 
	\overline{ \hat v(k, \eta')} 
	c_0^2 k^4 \eta_3 d k  d^2 \eta' \cr
	&=& (2 \pi)^3 \int_{-\infty}^\infty \int_{ |\eta'| < 1}
	 \hat u (k, \tilde \eta')  
	\overline{ \int_{ |\tilde \eta'|<1} 
	\overline{ \hat S(k, \eta', \tilde \eta') }
	 \hat v(k, \eta')  \eta_3  d^2 \eta' }  
	c_0^2 k^6 dk d^2 \tilde \eta'  \cr & & \label{adj1}
\end{eqnarray}
From (\ref{kernel}), (\ref{reciprocity}), and (\ref{symmetry}), we have 
\begin{eqnarray}
{( \mathcal P \mathcal S \mathcal P U, V )_W \over (2 \pi)^3}  &=&
	 \int_{-\infty}^\infty 
	\int_{|\eta'| <1} \hat u (k, \tilde \eta')  
	\overline{ \int_{ |\tilde \eta'| <1} 
	\overline{ - (i k / 2) A(k, -\tilde \eta^-, - \eta^+) }
	 \hat v(k, \eta') d^2 \eta' }c_0^2 k^6 d k d^2 \tilde \eta'\cr
 	&=& \int_{-\infty}^\infty 
	\int_{|\eta'|<1} \hat u (k, \tilde \eta')  
	\overline{ \int_{|\tilde \eta'|<1} (i k / 2) 
	A(-k, -\tilde \eta^-, - \eta^+) 
	 \hat v(k, \eta') d^2 \eta' }c_0^2 k^6 d k  d^2 \tilde \eta'\cr
	&=& \int_{-\infty}^\infty \int_{|\tilde \eta'|<1} 
	\hat u (k, \tilde \eta')  
	\overline{ \int_{|\eta'| <1}  
	\hat S(-k, -\tilde \eta', - \eta') 
	 \hat v(k, \eta') d^2 \eta' } \tilde \eta_3(-k) c_0^2 k^6
	 d k d^2  \tilde \eta' , \cr & & \label{adj2}
\end{eqnarray} 
where in the last equality we have used the fact that
$ \hat S(-k, \eta', \tilde \eta') = ik(2 \eta_3)^{-1} 
	A(-k, \eta^+, \tilde \eta^-)$.  
We note that in the course of this
computation, the $\eta_3$ has disappeared and has been replaced by
$\tilde \eta_3$; this is because of the $\eta_3$ in the denominator of
(\ref{kernel}).

We see that the action of the adjoint 
$(\mathcal P \mathcal S \mathcal P)^*$ on V is 
given in the transform domain by the $\eta'$ integral of (\ref{adj2}):
\begin{equation}
F {\mathcal F}^{-1} (\mathcal P \mathcal S \mathcal P)^* \mathcal F F^{-1}
	 \hat v ( k, \tilde \eta') =
	 \int_{| \eta'| <1} \hat S(-k, -\tilde \eta', - \eta') 
	 \hat v(k, \eta') k^2 d^2 \eta'. \label{adjS}
\end{equation}

However, $\hat S(k, \eta', \tilde \eta') $ is defined by 
\begin{equation}
\delta \left(c_0(k - \tilde k) \right) \hat S(k, \eta', \tilde \eta') =
	\int_{-\infty}^\infty \int_{-\infty}^\infty
	 \int \int e^{-i k \eta' \cdot x'} e^{i k c_0 t} 
	\mathcal S (t - \tau, x', y')
	e^{-i \tilde k c_0 \tau} e^{ik \tilde \eta' \cdot y'} d^2 x' 
	d^2 y' d \tau dt. \label{hatS}
\end{equation}
From this, we see that
\begin{equation}
\delta \left( c_0(\tilde k - k) \right) \hat S(-k, -\eta', -\tilde \eta') =
	\int_{-\infty}^\infty \int_{-\infty}^\infty	
	 \int \int e^{-i k \eta' \cdot x'} e^{-i k c_0 t} 
	\mathcal S (t - \tau, x', y')
	e^{i \tilde k c_0 \tau} e^{ik \tilde \eta' \cdot y'} d^2 x' 
	d^2 y' d \tau dt. \label{hatS-}
\end{equation}
Letting $t \rightarrow -t$ and $\tau \rightarrow -\tau$ in (\ref{hatS-})
shows that the kernel of (\ref{adjS}) corresponds to 
the operator $T \mathcal P \mathcal S \mathcal P T$, given by 
\begin{equation}
(T \mathcal P \mathcal S \mathcal P T)  V(t,x') 
	= \int_{-\infty}^\infty \int_{-\infty}^\infty 
	\mathcal S(\tau - t, x',y') V(t,y') 
	d\tau d^2 y'.
\end{equation}
Thus $(\mathcal P \mathcal S \mathcal P)^* 
= T (\mathcal P \mathcal S \mathcal P) T$.
\medskip

\subsection{Compactness}

\begin{theorem} Assume that the sound speed $c(x)$ is bounded and
differs from $c_0$ only in a 
bounded subset of the half-space $x_3 \leq -h < 0$.
Then the fixed-frequency scattering operator 
$\hat S$ is compact on the weighted space 
\begin{equation}
L^2_g = \{ f: f(\eta') |1 - |\eta'|^2|^{1/4} \in L^2 \}.
\end{equation}
\end{theorem}

Proof.  We use (\ref{kernel}) and (\ref{xi3}) to compute the square of the
Hilbert-Schmidt norm of $\hat S$ in the space $L^2_g$:
\begin{equation}
\| \hat S \|^2_{H.S.} = 
	{k^4 \over 4} \int \int 
	\left| \int e^{-i k \eta \cdot y} V(y) \psi (k, y, \tilde \eta') 
	d^3 y \right|^2
	\left| {1 - |\tilde \eta'|^2 \over 1 - | \eta'|^2} \right|^{1/2}
	d^2 \eta' d^2 \tilde \eta'.
	\label{HS}
\end{equation}

We use the fact that $\psi$ can be split into an incident and scattered field
via (\ref{LS}).  This allows us to split the kernel (\ref{kernel}) into two parts:
\begin{equation}
\hat S(k, \eta, \tilde \eta') = 
	\hat S_B(k, \eta, \tilde \eta') + 
	\hat S_{sc}(k, \eta, \tilde \eta')
\end{equation}
where the ``Born" term is
\begin{equation}
\hat S_B(k, \eta, \tilde \eta') =  -{i k\over 2 \eta_3} 
	\int e^{ik (\tilde \eta  - \eta) \cdot y} V(y) d^3 y \label{kernelB}
\end{equation}
and
\begin{equation}
\hat S_{sc}(k, \eta, \tilde \eta') =  -{i k^3 \over 2 \eta_3} 
	\int e^{-ik \eta \cdot y} V(y)  
	\int g(y-z) V(z) \psi (k, z, \tilde \eta') d^3 z d^3 y. 
	\label{kernelsc}
\end{equation}
We compute the Hilbert-Schmidt norm of each part.

The Hilbert-Schmidt norm of the Born term is
\begin{equation}
\| \hat S_B \|^2_{H.S.} = {k^2 \over 4} \int \int 
	\left| \int e^{ik(\tilde \eta - \eta) \cdot y} V(y)  d^3 y \right|^2
	\left| {1 - |\tilde \eta'|^2 \over 1 - | \eta'|^2} \right|^{1/2} 
	 d^2 \eta' d^2 \tilde \eta' . 	\label{HSB}
\end{equation}
The $y$ integral of (\ref{HSB}) can be written
\begin{equation}
\int e^{ik (\tilde \eta' - \eta') \cdot y'} 
	e^{-ik (\tilde \eta_3 + \eta_3) y_3} V(y) d^3y.
\end{equation}
For $|\eta'|<1$, for which $\eta_3$ is real, this integral is bounded
when $V$ is integrable.  For $|\eta'| > 1$, the integral decays
exponentially because $V$ is supported in the region where $y_3 \leq h
< 0$.  The same comments apply to the behavior in $\tilde \eta'$.  Thus
this integral is bounded by $c \exp (-hk ( |\eta'|  + |\tilde \eta'|))$.
This estimate can easily be used to show that $\| \hat S_B \|_{H.S.}$
is finite.

The Hilbert-Schmidt norm of the scattered part is
\begin{equation}
\| \hat S_{sc} \|^2_{H.S.} = {k^6 \over 4} \int \int 
	\left| \int e^{-ik \eta \cdot y} V(y) 
	\int g(y-z) V(z) \psi (k, z, \tilde \eta') d^3 z d^3 y \right|^2
	\left| {1 - |\tilde \eta'|^2 \over 1 - | \eta'|^2} \right|^{1/2} 
	 d^2 \eta' d^2 \tilde \eta' . 	\label{HSsc}
\end{equation}

An application of the Cauchy-Schwarz inequality shows that the $y$ integral of 
(\ref{HSsc}) is bounded by
\begin{equation}
\|  e^{-ik \eta \cdot y} V(y) |^2 \|_{L^2(y)} 
	\left\|  |V (y)|^{1/2} \int g(y-z) V(z) \psi(k, z, \eta')
	 d^3z \right\|_{L^2(y)}
	 \label{HSscest}
\end{equation}
Direct computation shows that the first norm in (\ref{HSscest}) is
bounded by $c \exp (-hk | \eta'|)$.  Standard scattering theory arguments
(See Appendix C) can be used to show that the second norm appearing in
(\ref{HSscest}) is bounded by $c \exp (- hk |\tilde \eta'| )$; thus the
$y$ integral satisfies the same bounds as in the Born term.  QED

{\bf Remark.}  In \cite{CI}, it was shown that on $L^2_g$, 
the operator $P\hat S P$ has norm less than or equal to one.
 
For future reference, we note that the operator 
$a(k) = ( P S P)^*( P S P)$ is given explicitly as
\begin{equation}
(F a f )(\tilde \eta') =  \int_{|\tilde \eta'| <1} \int_{|\zeta'| <1}  
	\hat S(-k, -\tilde \eta', - \zeta') 
	\hat S(k, \zeta',  \eta')  d^2 \zeta'
	 \hat f( \eta') k^4 d^2 \eta' . \label{a}
\end{equation}

\medskip

\section{Appendix:  Convergence of the Iterative Algorithm}

For simplicity of notation we consider only the case when $\mathcal S_0 = 0$.
In this case, the $n$th iterate is
\begin{equation}
U^\downarrow_n=
	{(T (\mathcal P \mathcal S \mathcal P) T 
	(\mathcal P \mathcal S \mathcal P))^n U^\downarrow_0 \over
	( (T (\mathcal P \mathcal S\mathcal P) T  
	(\mathcal P \mathcal S\mathcal P) )^n
	 U^\downarrow_0 , \Psi_B )_W }
	= {\mathcal A^n U^\downarrow_0 \over
	( \mathcal A^n U^\downarrow_0 , \Psi_B )_W},
	\label{Un}
\end{equation}
where  $\mathcal A= (\mathcal P \mathcal S \mathcal P)^* 
(\mathcal P \mathcal S \mathcal P)$, where the star denotes the
adjoint with respect to the flux inner product $(\cdot,\cdot)_W$.
We simplify the notation by dropping the arrow on $U$.  
Because we expect the limit to be a distribution,
we consider the quantity  
$ ( U_n, \Phi )_W 
= \int ( u_n, \phi )_{x'} c_0^2 k^4 dk$, 
where $\Phi$ is a smooth test function, $\phi =  {\mathcal F} \Phi$,
and $(\thinspace , )_{x'}$ denotes the weighted inner product 
that for smooth functions is 
$(f, \phi )_{x'} = (\hat f, \hat \phi )_{L^2_g}=
\int \hat f(\eta') \overline{\hat \phi(\eta')} |\eta_3| d^2\eta'$.

Specifically, we consider test function that are functions of
$t$ and $x'$.  When inverse Fourier transformed in $t$ and
Fourier transformed in space, at each frequency they must be
in $L^2_g$, and in the frequency variable they
must be integrable and (uniformly) supported in the
compact interval $[-B, B]$.  We denote
this space of test functions by $X$.

In the frequency domain, the $n$th iteration is 
${\mathcal F}^{-1}(\mathcal A^n U) = a^n u$, 
where 
$\mathcal F^{-1}$ denotes the inverse Fourier transform (\ref{tFT}) and 
$a = {\mathcal F}^{-1} \mathcal A {\mathcal F}= (PSP)^*(PSP)$. 
From our normalization of $U_n$, we have
\begin{equation}
( U_n, \Phi_B )_W = { (\mathcal A^n U, \Phi_B )_W \over 
	(\mathcal A^n U, \Psi_B)_W }
	= { \int \langle a^n u, \phi_B \rangle_{x'} c_0^2 k^4 dk
	\over \int \langle a^n u, \psi_B \rangle_{x'} c_0^2 k^4 dk }. 
	\label{UnPhi}
\end{equation} 

We note that $\mathcal A$ and $a(k)$ are self-adjoint on the space $W$
 and on $L^2_g$, respectively.  Moreover, Theorem 1
shows that $a(k)$ is compact on $L^2_g$, and it can therefore be
written $a(k) = \sum_l \lambda_l(k) P_l (k)$, where
$\lambda_l \geq \lambda_{l+1}$ and the $P_l$ are orthogonal projections.
Because $\mathcal A$ is non-negative, all the eigenvalues $\lambda_l$
are non-negative.  We see from (\ref{a}) that $a$ is analytic in $k$,
and the $\lambda$s are therefore piecewise analytic \cite{K}.

Suppose that $\lambda_0(k)$ attains its maximum in the set
$\{ |k| < B \}$ at $k_0$, and that $\lambda(k) = M$.  
Then in a neighborhood of $k_0$, $\lambda(k)$ has a 
Taylor expansion whose first two terms are $M - b(k - k_0)^p$
for some $b$ and some integer $p$.  We call $p$ the {\it order} 
of $\lambda_0$.

We allow eigenvalues with different indices to coincide at a point;
thus it is possible that a finite number of eigenvalues also attain
the maximum $M$ at $k_0$.  In this case, these eigenvalues
have a Taylor expansion similar to that of $\lambda_0$, possibly with
different $b$s and $p$s.  In this case we also refer to the relevant
$p$ as the order of the eigenfunction.

We will need the following lemma.

\begin{lemma}
Assume that $b$ is positive and that $p$ is an integer.  Then for large $n$, 
\begin{equation}
I(n,p) = \int_0^h (1 - b k^{p})^n dk \sim {C(p) \over  (bn)^{1/p}},
\end{equation}
where $C(p)$ is a nonzero constant independent of $n$.  
Thus the convergence to zero of $I(n,p)$ is slower for larger $p$
and smaller $b$.
\end{lemma}

Proof. 
Let $s = b^{1/p} k$.  Then $I = b^{-1/p} \int_0^{b^{1/p} h} (1 - s^p)^n ds$.
Replacing the upper limit by 1 results in an error that is exponentially small in $n$.  
Denote by $I_n$ the integral $\int_0^1 (1-s^p)^n ds$.  Then we can write
\begin{equation}
I_{n+1} = \int_0^1 (1-s^{p}) (1-s^{p})^n ds = I_n - \int_0^1 s^{p} (1-s^{p})^n ds. \label{L.1}
\end{equation}
In the integral of (\ref{L.1}), we integrate by parts, differentiating $s$ 
and integrating $s^{p-1} (1-s^{p})^n$.  
The boundary term vanishes, and (\ref{L.1}) becomes
\begin{equation}
I_{n+1} = I_n + {1 \over p(n+1)} I_{n+1}.
\end{equation}
Solving for $I_{n+1}$ gives the recursion 
\begin{equation}
I_{n+1} = { p(n+1) \over p(n+1) +1} I_n.
\end{equation}
Since $I_0 =1$, we have
\begin{equation}
I_n = \left( {p \over p+1} \right) \left( {2p \over 2p +1}\right) \cdots 
	\left( {np \over np+1} \right).
\end{equation}
Taking reciprocals and logs and expanding, we find that
\begin{eqnarray}
- \log I_n &=& \sum_{j=1}^n \log (1 + 1/(jp)) \cr
	&=& {1 \over p} \sum_{j=1}^n {1 \over j} - \sum_{m=2}^\infty 
		{(-1)^m \over m p^m} \sum_{j=1}^n {1 \over kj^m}
\end{eqnarray}
Exponentiating and taking reciprocals again, we have
\begin{equation}
I_n = C(n,p) \exp( - {1\over p} \sum_{j=1}^n {1 \over j}) , \label{L.2}
\end{equation}
where $C(n,p)$ has the large-n limit
\begin{equation}
\lim_{n \rightarrow \infty} C(n,p) = 
	\exp \left(- \sum_{m=2}^\infty { (-1)^m \over m p^m} \zeta(m) \right),
\end{equation}
where $\zeta(m)$ denotes the Riemann-zeta function 
$\zeta(m) = \sum_{k=1}^\infty k^{-m}$.  Thus we see that the large-n behavior
of $I_n$ is determined by the second factor of (\ref{L.2}).  

We determine the large-n behavior of this second factor as follows.
From approximating the sum by a Riemann integral, we have the estimate
\begin{equation}
\log n \leq \log (n+1) \leq \sum_{j=1}^n {1\over j} \leq 1 + \log n
\end{equation}
We multiply by $-1/p$ and exponentiate to obtain 
\begin{equation}
e^{-1/p} n^{-1/p} \leq \exp( - {1\over p} \sum_{j=1}^n {1 \over j}) 
	\leq n^{-1/p}
\end{equation}
QED 

\begin{theorem}
Assume that $a$ is an analytic self-adjoint-compact-operator-valued function of $k$
having the representation 
$a(k) = \sum_l \lambda_l(k) P_l(k)$, 
where $\lambda_l \geq \lambda_{l+1}$
and the $P_l$ are orthogonal projections.   
Assume that $\lambda_0$ is not a constant function of $k$.  
Then for test functions  $\phi_B(k, x')$ and  
$\psi_B(k, x')$ in $X$ whose support in the frequency domain
is in the set $\{ |k| \leq B \}$,
\begin{equation}
\lim_{n \rightarrow \infty}
	{\int (a^n  u, \phi_B )_{x'} k dk \over
	 \int (a^n u, \psi_B )_{x'} k dk} 
	= {\sum_{l,j} \beta_{lj} ( P_l u, \phi_B )_{x'} 
		(k_j) \over
	\sum_{l,j}\beta_{lj} (P_l u, \psi_B )_{x'} (k_j)}
\end{equation}
where the sums are over those indices $j$ and $l$ for which
$\lambda_l(k_j) = M$, where $M$ is the maximum of $\lambda_0$ in the set
$\{ k : |k| \leq B\}$, and for which
$\lambda_l$ has maximal order at $k_j$.  
\end{theorem}

\medskip
\leftline{Proof}
The representation for $a$ allows us to write (\ref{UnPhi}) as 
\begin{equation}
(U_n, \Phi_B )_W 
	= {\int  \sum_l \lambda^n_l (P_l u, \phi_B )_{x'} 
		k dk \over
	\int \sum_l \lambda^n_l (P_l u, \psi_B )_{x'} k dk}
		\label{U_n}
\end{equation}
The $\lambda_l$ and $P_l$ are peicewise analytic functions of $k$ \cite{K}.
In particular $\lambda_0(k)$ is piecewise analytic, and therefore attains
its maximum $M$ on a discrete subset of the set $\{ |k| \leq B \}$.  
We cover the support of $\phi$ with open intervals $N_j$ so that
each $N_j$ contains only one $k_j$.  We decompose the
test function $\phi_B$ as $\phi_B = \sum_j \phi_j$ \cite{H}, 
where the $\phi_j$ are in $C^\infty_0(N_j)$ and $\phi_j = \phi_B$ in a
neighborhood of $k_j$.  We carry out a similar decomposition for
$\psi_B$.  

With the notation $f_{l,j}(k) = k ( P_l u, \phi_j )_{x'}$ and
$g_{l,j}(k) = k ( P_l u, \psi_j )_{x'}$, we can write (\ref{U_n}) as
\begin{equation}
( U_n, \Phi_B )_W 
	= { \sum_{l,j} \int_{N_j}\lambda^n_l (k) f_{l,j}(k) 
		 dk \over
	 \sum_{l,j}\int_{N_j} \lambda^n_l (k) g_{l,j} (k) 
		 dk}.  \label{B.1}
\end{equation}

We divide the numerator and denominator of (\ref{B.1}) by $M^n$, and write 
$r_l(k) = \lambda_l (k) / M$; thus  
$|r_0| \leq 1$, and $|r_l| < 1$ for all but a finite number of
values of $l$. 
Then (\ref{B.1}) can be written
\begin{equation}
( U_n, \Phi_B )_W =
	{\sum_{l,j} \int_{N_j} r_l^n (k) f_{l,j} (k) dk \over
	\sum_{l,j} \int_{N_j} r_l^n (k) g_{l,j} (k) dk }. 
	\label{quotient}
\end{equation}

We write
\begin{equation}
I^n_{l,j} = \int_{N_j} r_l^n (k) f_{l,j} (k) dk. \label{Ilj}
\end{equation}

We multiply and divide $I^n_{l,j}$ by
$\int_{N_j} r^n_0 (k) dk$, and 
write $\zeta^n_l = r^n_l / \int_{N_j} r^n_0 dk$.

For those $l= 0, 1, \ldots L_j $ for which $\lambda_l$ 
attains the maximum $M$ and
thus $r_l$ attains the value $1$, we 
add and subtract $f_{l,j}(k_j)$ to the quotient, obtaining
\begin{equation}
I_{l,j}^n = \left( f_{l,j} (k_j) + \int_{N_j} \zeta^n_l(k)  
	(f_{l,j}(k) - f_{l,j} (k_j)) dk \right) 
	\int_{N_j} r_0^n (k) d k  \label{I_0}
\end{equation}
We will show that the integral term within the 
parentheses on the right side of (\ref{I_0}) 
vanishes as $n$ goes to infinity.  
For this we use the two facts that 1) except at
$k = k_j$, $\zeta^n_l $ converges to zero pointwise as $n$ goes to 
infinity; and 2) $\int \zeta^n_l = 1$ for all $n$.
  
Given $\epsilon>0$, we choose $N^\epsilon_j$ so small that on $N^\epsilon_j$,
$|f_{l,j}(k) - f_{l,j}(k_j)| < \epsilon/2$. 
The integral in parentheses on the right side of (\ref{I_0})
we split into two integrals, namely $A^n$ and $B^n$, where 
\begin{equation}
A^n = \int_{N^\epsilon_j} \zeta^n_l (k)
	(f_{l,j} (k) - f_{l,j}(k_j)) dk,
\end{equation}
and 
\begin{equation}
B^n = \int_{N_j \setminus N^\epsilon_j}\zeta^n_l (k)
	(f_{l,j} (k) - f_{l,j}(k_j)) dk.
\end{equation}
Then $A^n < \epsilon/2$.  Next, we choose $N$ so large that for $n$
greater than $N$, $B^n < \epsilon/2$.  
This shows that the integral in parentheses on the right side of 
(\ref{I_0}) vanishes as $n$ goes to infinity.  Thus (\ref{I_0}) is a product
of a factor converging to  $f_{l,j}(k_j)$ and a factor converging to
zero.

For $l= L_j+1, L_j+2,  \ldots $, for which $\lambda_l$ 
is strictly less than $M$, we write
\begin{equation}
I^n_{l,j} = \left( \int \zeta^n_l (k) f_{l,j} (k) dk \right)
	\int_{N_j} r_0^n (k) dk.
\end{equation}
In this case, $\zeta^n_l$
converges to zero pointwise for all $k$.  Thus the integral in 
parentheses of $I^n_{l,j}$ converges
to zero by the Lebesgue Dominated Convergence Theorem. 

To estimate the tail of the sum over $l$ in (\ref{quotient}), 
we choose $l_0$ so large
that for $l>l_0$, $\zeta_l(k) < 1/2$ for all $k$ in $N_j$.  
This is possible because
the the compactness of $a$ implies that its eigenvalues decrease to zero.
Thus for each $j$ we have 
\begin{equation}
\left| \sum_{l>l_0} I^n_{l,j} \right| 
	= \left| \sum_{l>l_0} \int_{N_j} \zeta_l^n (k)  
	( P_l u, \phi_j )_{x'} k dk  \right|
	\leq {1 \over 2^n} \int_{N_j} \sum_{l>l_0} \left| (u, \varphi_l)_{x'} 
	(\varphi_l, \phi_j)_{x'} k  \right| dk ,
\end{equation}
where the $\varphi_l$ are the normalized eigenfunctions of $a$.  Here it
may be necessary to reindex the sum.  Each
of the sequences $(u, \varphi_l)_{x'}$ and $(\varphi_l, \phi_j)_{x'}$
is in $l^2$, and the inner product of two $l^2$ sequences is in $l^1$.
We therefore find that the sum over $l$ is bounded by 
$\| u \|_k  \| \phi_j \|_k$, where $\| \cdot \|_k$ denotes
the norm in the space $L^2_k$.  Thus we have
\begin{equation}
\left| \sum_{l>l_0} I^n_{l,j} \right|
	\leq {1 \over 2^n} \int_{N_j} \| u \|_k  \| \phi_j \|_k
	k dk
\end{equation}
which shows that the tail of the sum converges
to zero as $n$ goes to infinity. 

The same arguments, of course, apply to the denominator of (\ref{quotient}). 
Thus we see that the leading order behavior of (\ref{quotient}) is
given by the expression
\begin{equation}
( U_n, \Phi_B )_W \sim
	{ \sum_{l,j} f_{l,j} (k_j)  \int_{N_j} r_0^n (k) d k  
	\over \sum_{l,j} g_{l,j}(k_j) \int_{N_j} r_0^n (k) dk },
\end{equation}
where the sum in $l$ is over those values for which $\lambda_l$
attains the maximum $M$ at $k_j$.
The terms $\int_{N_j} r_0^n (k) d k $, however, go to
zero for large $n$.  We must therefore consider their behavior in
more detail.

As we have seen, in the neighborhood of
$k = k_j$, $r_0(k)$ has an expansion of the form
$r_0(k) = 1 - b_j(k - k_j)^{p_j} + \ldots$, where the
positive integer $p_j$ is the order of the eigenfunction at 
$k_j$.  The lemma shows that the order $p_j$ controls the speed with which
$\int_{N_j} r_0^n (k) dk$ goes to zero with $n$:  the larger $p_j$,
the more slowly $\int r_0^N dk$ coverges to zero.  
We divide the numerator and denominator  of (\ref{quotient}) by 
$\int_{N_j} r_0^n dk$ corresponding to the slowest decay. 
Finally, we take the limit of the resulting quotient as $n$ goes to infinity.  
This shows that the
quotient (\ref{quotient}) converges to
\begin{equation}
( U_n, \Phi_B )_W \sim
	{\sum_{l,j} \beta_j (P_l u, \phi )_{x'} (k_j) \over
	\sum_{l,j} \beta_j (P_l u, \psi )_{x'} (k_j)}
	= \int \int \sum_{l,j} \beta_j {P_l u (k_j, x') \over
	\sum_{l,i} \beta_i (P_l u, \psi )_{x'} (k_i)} 
	\delta_{k_j} (k) \phi (k, x') d^2 x' dk, 
\end{equation}
where the $\beta_j$ are proportional to 
$\beta_j = k_j C(p_j) / b_j^{1/p_j}$,
and where the sums are over those indices $j$ and $l$ for which $\lambda_l$
has maximal order at $k_j$.

QED.

\begin{corollary} 
Assume that the sound speed $c(x)$ is bounded and
differs from $c_0$ only in a 
bounded subset of the half-space $x_3 \leq -h < 0$.  Then
for test functions $\Psi_B,\Phi_B$ in $X$ whose 
Fourier transforms with respect to time are supported in
$-B \leq k \leq B$, 
$ U_n $ as defined by (\ref{Un}) converges to 
\begin{equation}
{1 \over 2 \pi \sum_{l,i} \beta_k (P_l u, \psi_B )_{x'} (k_i)} 
\sum_{l,j} \beta_j P_l u (k_j, x') e^{- i k_j c_0 t},
\label{tresult}
\end{equation}
where the sums  are over those indices $j$ and $l$ for which 
$\lambda_l(k_j)$ attains the maximum $M$ and has maximal
order.  
\end{corollary}

Proof.  To apply Theorem 2, we need only check that 
$a = \mathcal F^{-1} (\mathcal P \mathcal S \mathcal P)^* 
(\mathcal P \mathcal S \mathcal P) \mathcal F$ is
indeed an analytic compact-operator-valued function of $k$.
Analyticity was shown in section A.1; compactness was shown in
section A.3.   The largest eigenvalue $\lambda_0$ cannot be constant:
$S$ is zero at $k =0$, which implies that all the $\lambda_j$
are zero there.  Thus if $\lambda_0$ were constant it would be zero,
and $S$ itself would be zero.

We note that both the numerator and denominator in (\ref{tresult}) 
can be zero, in which case (\ref{tresult}) is not defined. However
we do not study this case since the denominator is non-zero
for a generic test function $\Psi_B$.
QED

\section{Appendix:  Results from ``Standard" Scattering Theory}

The solution of the equation 
\begin{equation}
(\nabla^2 + k^2 - V(x)) \psi(k, x) = 0
\end{equation}
corresponding to an incident plane wave and a scattered field satisfying
outgoing boundary conditions satisfies the Lippmann-Schwinger 
integral equation
\begin{equation}
\psi(k, x, \tilde \eta') = \exp (ik \tilde \eta \cdot x) - 
	k^2 \int g(k,x,y) V(y) \psi(k, y, \tilde \eta') d^3 y, 
	\label{C.LS}
\end{equation}
where $g$ is the usual outgoing Green's function (\ref{Green}).

The initial difficulty with the Lippmann-Schwinger equation is that the 
incident field
has infinite energy.  This difficulty, however, can be circumvented by
 multiplying the
whole equation by $|V(x)|^{1/2}$ \cite{S}.  This converts (\ref{C.LS}) into
\begin{equation}
\zeta(k, x, \tilde \eta') = \zeta_0(k, x, \tilde \eta') + 
	k^2 \int K(k, x-y) \zeta (k, y, \tilde \eta') d^3y, 
	\label{factoredLS}
\end{equation}
where
\begin{eqnarray}
\zeta(k, x, \tilde \eta') &=& 
	|V(x)|^{1/2} \psi (k, x, \tilde \eta'), \\
\zeta_0(k, x, \tilde \eta') &=& 
	|V(x)|^{1/2} e^{ik \tilde \eta \cdot x}, \\
V_{1/2}(y) &=& V(y) / |V(y)|^{1/2} \\
K(x,y) &=& |V(x)|^{1/2} g(k, x,y) V_{1/2}(y) \\
\end{eqnarray}

When $V$ has compact support, (\ref{factoredLS}) is an integral equation on a
bounded region.  The kernel $K$ is an Hilbert-Schmidt-valued function that
is analytic in the entire complex $k$-plane.
By the analytic Fredholm theorem \cite{RS1}, the integral equation (\ref{factoredLS}) is
therefore uniquely solvable everywhere except at a discrete set of values of $k$,
(the ``exceptional points'')
and moreover the solution $\zeta$ is a meromorphic function of $k$ with poles at
these exceptional points.  In addition, the arguments of   
 \cite{A} \cite{RS4} show that the only possible real exceptional point is $k=0$.
However, when $k=0$, we also have $k=0$, and (\ref{factoredLS}) reduces 
to the equation $\zeta = |V|^{1/2}$.  Thus the operator $(I-k^2 K)^{-1}$ is analytic
in a neighborhood of the real $k$-axis.  

This argument shows that for each $k$, $\zeta = |V|^{1/2} \psi$ is in $L^2$.  
Then the quantity needed
in section A.3, namely $\|  |V|^{1/2} \int gV \psi \|$, can be rewritten as 
$\| K \zeta \| =\|  K (I-K)^{-1} \zeta_0 \|  \leq \| K \| \|(I-K)^{-1}\| \| \zeta_0 \|$.  
Moreover, by explicit computation, we see that 
$\| \zeta _0(k, \cdot, \tilde \eta') \| \leq c e^{-hk |\tilde \eta'| }$, 
where $V$ is supported in the region $ y_3 \leq -h < 0$.

\end{document}